\documentclass{amsart}
\usepackage{marktext}
\usepackage{gimac}

\newtheorem{question}[theorem]{Question}
\DeclareMathOperator{\cov}{cov}
\newcommand{\loc}{\mathrm{loc}}

\renewcommand{\curl}{\nabla^\perp \cdot}
\DeclareMathOperator{\crl}{curl}
\DeclareMathOperator{\divv}{div}

\begin{document}
% Title etc.
\title[Euler equations with random initial vorticity]{Incompressible 2D Euler equations with non-decaying random initial vorticity}
\begin{abstract}
  % Tile the plane $\R^2$ with squares of sides of unit length.
  % In each square toss a fair coin and assign constant vorticity of strength $+1$ if heads, $-1$ if tails.
  % In this paper we prove global well-posedness of weak solutions of the incompressible Euler equations for almost every initial vorticity described above.
  % The main contribution of our work is the construction of a
  % %unique, up to a constant, divergence-free vector field in the
  % corresponding initial velocity field that grows slowly at infinity, 
  % %plane having this given vorticity,
  % which enables us to apply  a recent well-posedness result of Cobb and Koch.
  Consider a random initial vorticity~$\omega_0(x) = \sum_{n\in \mathbb{Z}^2} a_n \phi(x-n)$, where $\phi$ is bounded and compactly supported and $\{a_n\}$ are independent, uniformly bounded, mean~$0$, variance~$1$ random variables (i.e.\ $\omega_0$ is an array of randomly weighted vortex blobs).
  We prove global well-posedness of weak solutions to the Euler equations in~$\R^2$ for almost every such initial vorticity.
  %that for almost every such initial vorticity, the unforced  global well-posedness of weak solution to the unforced incompressible two-dimensional Euler equations in~$\R^2$ with 
  %For almost all realizations of the $\{a_k\}$, we construct an initial  velocity field $u_0$ with divergence zero, curl equal to the corresponding $\omega_0$, $u_0(0) = 0$ and which satisfies a suitable integral growth estimate. We then use a recent result by D. Cobb and R. Koch \cite{CobbKoch24} for existence and uniqueness of a weak solution with the realization $u_0$ as initial velocity.
  The main contribution of our work is the construction of a
  %unique, up to a constant, divergence-free vector field in the
  corresponding initial velocity field that grows slowly at infinity, 
  %plane having this given vorticity,
  which enables us to apply  a recent well-posedness result of Cobb and Koch.
\end{abstract}
\author[Iyer]{Gautam Iyer}
\address{%
  Department of Mathematical Sciences, Carnegie Mellon University, Pittsburgh, PA 15213.
}
\email{gautam@math.cmu.edu}

\author[Lopes Filho]{Milton C. Lopes Filho}
\address{Instituto de Matem\'atica, Universidade Federal do Rio de Janeiro, Cidade Universit\'aria -- Ilha do Fund\~ao, Caixa Postal 68530, 21941-909 Rio de Janeiro, RJ -- BRAZIL}
\email{mlopes@im.ufrj.br}

\author[Nussenzveig Lopes]{Helena J. Nussenzveig Lopes}
\address{Instituto de Matem\'atica, Universidade Federal do Rio de Janeiro, Cidade Universit\'aria -- Ilha do Fund\~ao, Caixa Postal 68530, 21941-909 Rio de Janeiro, RJ -- BRAZIL}
\email{hlopes@im.ufrj.br}

\thanks{This work has been partially supported by
  the National Science Foundation through grants DMS-2342349, DMS-2406853, 
  CNPq through grants 304990/2022-1, 305309/2022-6,
  FAPERJ, through  grants E-26/201.209/2021, E-26/201.027/2022,
  and the Center for Nonlinear Analysis.
}
\subjclass{%
  Primary:
    76B47. %	Vortex flows for incompressible inviscid fluids
  Secondary:
    76B03, %	Existence, uniqueness, and regularity theory for incompressible inviscid fluids
    76F99. %	Turbulence / None of the above, but in this section
  }
\keywords{Euler equations, vorticity, turbulence}
\maketitle

%\begin{todo}[Choose a better title]
%  Suggestions:
%  \begin{enumerate}
%    \item Solutions to the Euler equations with random vortex patch initial data.
%    \item The Euler equations with random vortex patches.
%    \item Existence of \dots
%  \end{enumerate}
%\end{todo}

\section{Introduction}\label{s:intro}

The Euler equations govern the evolution  of the velocity field of an ideal, incompressible fluid, and are given by
\begin{equation}\label{e:euler}
\left\{
\begin{aligned}
  \partial_t u + u \cdot \grad u + \grad p = 0 & \quad \text{ in } (0,+\infty)\times \R^2,\\
  \dv u = 0 &\quad \text{ in } [0,+\infty)\times \R^2
  \,.
\end{aligned}
\right.  
\end{equation}
Here~$u$ represents the fluid velocity, and~$p$ is the pressure.
%The purpose of this paper is to study the unforced incompressible 2D Euler equations in~$\R^2$ with bounded ``random'' initial vorticity.
Our aim is to study these equations when the initial vorticity~$\omega_0 = \crl u_0$ is an array of vortex blobs with random weights.
Explicitly, we suppose
\begin{equation}\label{e:omega0intro}
  \omega_0(x) = \sum_{n \in \Z^2} a_n \phi(x - n)
\end{equation}
where~$\phi$ is a bounded compactly supported function, and~$\set{a_n}$ are independent, uniformly bounded, mean~$0$, variance~$1$ random variables.
The main result (Theorem~\ref{t:existence}, below) shows global well-posedness of the Euler equations for almost every such initial vorticity.
%normalized to have mean~$0$ and variance~$1$.

The motivation for studying this problem stems from the desire to understand the generic interaction of vortices, a long standing problem in the theory of turbulence.
A well studied phenomenon is the \emph{energy cascade}: three dimensional fluid flows tend to mix by filamentation and transfer energy to higher frequencies.
High frequencies are rapidly dissipated by the viscosity and the balance between these phenomenon leads to a characteristic power-law energy spectrum~\cite{Frisch95,Obukhov49,Corrsin51,Kolmogorov41}.

In two dimensions, the energy cascade is qualitatively different.
Fluids tend to transfer energy to lower frequencies leading to a phenomenon known as the \emph{inverse cascade} \cite{KraichnanMontgomery80,Chorin94,Wayne11,GallayWayne05}.
Physically, one manifestation of the inverse cascade is \emph{vortex coalescence} -- the tendency of vortices to coalesce into a smaller number of larger vortices.
One of the main goals in the study of two-dimensional turbulent flows has been to understand and explain this inverse cascade and its connection with vortex coalescence.
%To quote directly from \cite{Wayne11}:
%\begin{quote}
%  One of the main goals in the study of two-dimensional turbulent flows has been to understand and explain this inverse cascade and in particular to explain the tendency of the vorticity to coalesce into a smaller and smaller number of larger and larger vortices.
%\end{quote}

The problem studied in this paper is motivated by the above phenomenon.
Several turbulence theories capture generic fluid behavior by combining the evolution equations with statistical models, see for instance \cite{Kupiainen11} and references therein.
In this spirit, we study the ``generic'' interaction of vortices by starting with an array of vortex blobs with random weights, and evolving this configuration by the (unforced) incompressible Euler equations~\eqref{e:euler}.
We note that, if the $a_n$'s are i.i.d., then the initial vorticity ~\eqref{e:omega0intro} is statistically invariant under the action of $\Z^2$.
This is a step towards understanding spatial white-noise initial vorticity, which is a physically meaningful problem.
(We also mention that our result is different in flavor from the results of~\cite{AlbeverioCruzeiro90,Flandoli18,CruzeiroSymeonides18}, where the authors study solutions with low regularity and construct invariant measures.)

Presently, the main physical implication of our result (Theorem~\ref{t:existence}, below) is about velocity growth.
Reconstructing the initial velocity field from the vorticity~$\omega_0$ involves understanding the (linear) interaction of the vortex blobs.
The slow decay of the Biot--Savart kernel allows for long range interaction of the vortex blobs and may cause linear growth of the initial velocity field.
However, since the weights~$\set{a_k}$ are mean~$0$ and independent, there is a lot of cancellation in the interaction.
We use this to construct an initial velocity field that grows \emph{slower} than any positive power (in a certain integrated sense).
A recent result of Cobb and Koch~\cite{CobbKoch24} shows that this growth condition is preserved by the flow, and that~\eqref{e:euler} is globally well-posed.
We remark that the cancellation in the interaction in fact shows that the initial velocity field grows slower than the square-root of a logarithm on average (see Remark~\ref{r:optimality}, below), but this is not guaranteed to be preserved by the flow.

To explain the significance of this, and highlight the main mathematical difficulties involved, we briefly recall a few relevant results concerning the well-posedness of solutions to~\eqref{e:euler} for bounded vorticity.
When the initial vorticity is both bounded and integrable, the velocity is \emph{a priori} bounded and the well-posedness of~\eqref{e:euler} is now classical~\cite{Majda86,Judovic63}.
There are, however, a number of physical reasons to study~\eqref{e:euler} when the vorticity is merely bounded, and allow the velocity field to be unbounded.
This arises, for instance, in natural examples (e.g.\ Couette flow, which have a constant vorticity and a linearly growing velocity), or rigid rotations / hurricanes, where the velocity field grows linearly with the distance from the center, and in the context of oceanic flows far away from the coastline.

The study of merely bounded vorticity flows also poses a number of interesting mathematical challenges.
Serfati~\cite{Serfati95} (see also \cite{AmbroseKelliherEA15}) studied~\eqref{e:euler} when the initial velocity and the initial vorticity are both bounded.
Brunelli~\cite{Brunelli10} studied the case where the velocity may grow like
a square-root,
%~$\abs{x}^{1/2}$ as~$x \to \infty$,
but required an integral decay condition on the vorticity.
Without a decay assumption on the vorticity,
the first result allowing growth on the velocity was by Cozzi~\cite{Cozzi15}.
The logarithmic growth allowed in~\cite{Cozzi15} was then improved to a nearly square-root growth by Cozzi and Kelliher in~\cite{CozziKelliher19}.
%Cozzi~\cite{Cozzi15} allows a logarithmic growth on the velocity, and then 
%allowing logarithmic growth by E. Cozzi, see \cite{Cozzi15},
%conditions on the growth of the velocity were first relaxed
%then nearly square-root growth by E. Cozzi and J. Kelliher, see \cite{CozziKelliher19}.
%For reasons we will make clear, these results are not suited for our purpose.
In a deep and innovative recent work~\cite{CobbKoch24},  Cobb and Koch prove a global well-posedness result in a Morrey-type space which allows for nearly-square root growth of the initial velocity.
To the best of our knowledge, even local well-posedness of~\eqref{e:euler} with bounded vorticity, and allowing the velocity field to grow like a square-root, or faster, is open.
%obtained a well-posedness result where initial vorticity is bounded and the initial velocity satisfies a Morrey-type estimate, which includes near square-root growth, see \cite{CobbKoch24}. Because this condition is expressed in terms of an integral, the new result can be used to complete our analysis. 

In the context of our work, given the initial vorticity~\eqref{e:omega0intro}, consider a (random) incompressible velocity field $u_0$ for which~$\crl u_0 = \omega_0$.
Such a velocity field is unique up to the addition of the gradient of a harmonic function.
This velocity field may grow at infinity, and if the growth is like a square-root or faster, then well-posedness of~\eqref{e:euler} is unknown.

The growth of the velocity field is governed by the (linear) interaction of the vortex blobs in~\eqref{e:omega0intro}.
For simplicity, suppose~$\phi$ is the indicator function of the unit square~$[0, 1)^2$, and suppose~$\set{a_k}$ are~$\pm 1$ valued random variables.
If all the signs align on a large enough region (like a sector), it will force the velocity field~$u_0$ to grow linearly at infinity.
On the other hand, if the signs change frequently, there may be a lot of cancellation, and one may be able to find an initial velocity field~$u_0$ which grows very slowly, or even remains bounded.
In this paper we settle this question by producing a velocity field whose growth at infinity is (almost surely) slower than any positive power-law growth (in a Morrey-type sense).
The Cobb--Koch result~\cite{CobbKoch24} will now imply both global well-posedness of~\eqref{e:euler}, and that the slow initial growth of the velocity field is preserved for all time by the non-linear evolution of~\eqref{e:euler}.

\subsection*{Plan of this paper.}
The remainder of this article is divided into three sections. In Section \ref{mainres} we introduce some notation and state our main result. In Section \ref{proof} we present its proof.
In Section \ref{decorscaling} we obtain logarithmic bounds for the decorrelation scale of the initial velocity, we observe that the solutions we obtain with initial vorticity of the form $\omega_0(x/\epsilon)$ vanish as $\epsilon \to 0$ and we formulate related open problems.

\subsection*{Acknowledgements}
The authors wish to thank
  Nathan Glatt-Holtz, 
  Daniel Rodriguez Marroquin
  and
  Vlad Vicol
for helpful discussions.

\section{Main Result.} \label{mainres}

We begin by stating our main result.
For simplicity, we consider initial vorticity of the form
\begin{equation}\label{e:omega0def}
  \omega_0 \defeq \sum_{n \in \Z^2} a_n \one_{Q_n}
  \,,
\end{equation}
where~$Q_0 \defeq [0, 1)^2$ is the unit square, $Q_n = n + Q_0$, and~$\set{a_n}$ are random weights.

%Following the statement, we will define the spaces and the notion of weak solution which we require.
%Denote $\grad^\perp = (-\partial_2, \partial_1)$ and recall that the two-dimensional curl of a planar vector field $v$ is given by~$\grad^\perp \cdot v$. 

\begin{theorem}\label{t:existence}
  Suppose~$\set{a_n}$ are uniformly bounded, independent, random variables all of which have mean-zero and variance~$1$.
  For almost every initial vorticity $\omega_0$ defined by~\eqref{e:omega0def}, there exists a (random) velocity field~$u_0 \in L^2_\loc(\R^2)$ such that the following hold.
  \begin{enumerate}
    \item\label{i:divcurl}
      The distributional derivative of~$u_0$ satisfies~$\divv u_0 = 0$ and~$\crl u_0 = \omega_0$.
         %\item
    %  For every~$\alpha > 0$, we have~$u_0 \in L^2_\alpha$.
    %  Moreover, such a velocity field is unique in~$L^2_\alpha$ up to addition of a constant vector field.
    % Such an $u_0$ is unique up to the addition of a (random) vector field that is constant in space.

    \item
     There exists a unique weak solution of the 2D Euler equations~\eqref{e:euler} with  $u(0)=u_0$ such that for every~$\alpha \in (0, 1/2)$ we have
      \begin{equation}\label{e:uSpace}
	u \in L^\infty_\loc( [0, \infty);  L^2_\alpha(\R^2) )
	\,,
	\qquad
	\omega \in L^\infty( [0, \infty); L^\infty(\R^2))
	\,,
      \end{equation}
      and, also,  the far field condition,
      % Helena, you think there might be too many commas?
      \begin{equation} \label{e:farfieldcond}
	u(t) - u(0) \in \mathcal S'_h
	\,.
      \end{equation}
  \end{enumerate}

\end{theorem}

%\begin{remark}[Continuity of the velocity field]

%\end{remark}

The space~$L^2_\alpha$ appearing in~\eqref{e:uSpace} is a \emph{non-homogeneous, local Morrey space}, in the terminology of~\cite{CobbKoch24}.
In contrast to the traditional Morrey spaces~\cite{Stein93}, the space~$L^2_\alpha$ captures growth at infinity, and not local regularity.
Precisely, for $\alpha \geq 0$, the norm in~$L^2_\alpha(\R^2)$ is defined by
\begin{equation} \label{e:L2alphanorm}
  \norm{f}_{L^2_\alpha(\R^d)}
    % \defeq \sup_{R \geq 1} \frac{1}{R^{\alpha + d/p}} \norm{f}_{L^p(B_R)}
    % = \sup_{R \geq 1} \paren[\bigg]{\frac{1}{R^{\alpha p + d}} \int_{B_R} \abs{f}^p \, dx }^{1/p}
    \defeq \sup_{R \geq 1} \frac{1}{R^{\alpha + 1}} \norm{f}_{L^2(B_R)}
%= \sup_{R \geq 1} \frac{1}{R^{\alpha +1}} \paren[\bigg]{\int_{B_R} \abs{f}^2 \, dx }^{1/2}
    \,,
\end{equation}
where~$B_R \defeq \set{ x \in \R^2 \st \abs{x} < R}$.
Clearly any function for which~$\abs{f(x)} \leq C (1 + \abs{x}^\alpha)$ belongs to~$L^2_\alpha(\R^2)$.
(Here, and throughout this paper, we use~$C$ to denote a constant that may change from line to line.)

The space $\mathcal S'_h$ in \eqref{e:farfieldcond} is the Chemin space of homogeneous tempered distributions (see~\cite{Chemin04}) defined as follows.
Fix $\chi \in C^\infty_c(\R^2)$ such that $\chi \equiv 1$ in $B_1$, $0\leq \chi \leq 1$, and assume that the support of $\chi$ is contained in the ball $B_2$. Then:
\begin{equation} \label{e:Sprimeh}
\mathcal S'_h = \{ f \in \mathcal S' \st \chi(\lambda D) f\xrightarrow[]{\lambda \to \infty} 0  \text{ in } \mathcal S'\}.    
\end{equation}
Above, $D$ is the multiplier operator with symbol $i \xi$.

Lastly, recall from~\cite{CobbKoch24} the notion of weak solution used in Theorem~\ref{t:existence}.

\begin{definition}\label{d:weaksol}
Let $\alpha \geq 0$ and $T>0$.
  We say $u\in L^2_{\loc}([0,T);L^2_\alpha(\R^2))$ is a weak solution of the 2D Euler equations~\eqref{e:euler} with initial velocity $u_0 \in L^2_\alpha (\R^2)$ if the following hold.
\begin{enumerate}
    \item We have $\divv u = 0$ in $\mathcal D' ((0,T)\times\R^2)$.
    \item For every $\varPhi \in [C^\infty_c([0,T)\times \R^2)]^2$ such that $\divv \varPhi =0$ it holds that
      \[\int_0^T\int_{\R^2} \paren[\big]{u\cdot \partial_t \varPhi + u\cdot[(u\cdot\nabla)\varPhi]} \, dx \, dt   + \int_{\R^2} u_0(x) \cdot \varPhi(0,x) \, dx = 0.\]
\end{enumerate}
\end{definition}

\begin{remark}[Uniqueness of~$u_0$]
 By standard elliptic regularity, any velocity field~$u_0$ satisfying condition~(\ref{i:divcurl})  in Theorem~\ref{t:existence}, with bounded vorticity, will be continuous. If we additionally fix~$u_0$ at a point, and require~$u_0 \in L^2_\alpha$ for some~$\alpha < 1$, then this~$u_0$ is unique.
\end{remark}

\begin{remark}[Growth of~$u_0$]
  We will, in fact, construct~$u_0$  so that it grows like the square-root of a logarithm, on average.
  Explicitly Proposition~\ref{p:u0exist}, below, will imply that there exists a constant~$C > 0$ such that for every~$x \in \R^2$ we have
  \begin{equation}
    \E \abs{u_0(x)}^2 \leq C \ln( e + \abs{x})
    \,.
  \end{equation}
  This growth bound is expected to be optimal for reasons we elaborate on in Remark~\ref{r:optimality}, below, but does \emph{not} imply an almost sure bound of the form~$\abs{u_0(x)}^2 \leq \widetilde{C} \ln (e + \abs{x})$ for some (random) constant~$\widetilde{C}$.
\end{remark}
\begin{remark}[More general initial vorticity]
  The proof of Theorem~\ref{t:existence} applies in the following more general situation.
  Suppose~$\mathcal{D} \subseteq \R^2$ is a compact domain and~$\set{\phi_n \st n \in \Z^2}$ is a family of uniformly bounded functions that are all supported in~$\mathcal{D}$.
  As before let $\set{a_n \st n \in \Z^2}$ be a family of independent, uniformly bounded, mean~$0$ variance~$1$ random variables.
  If instead of~\eqref{e:omega0def} we define~$\omega_0$ by
  \begin{equation}\label{e:xi0general}
    \omega_0(x) \defeq \sum_{n \in \Z^2} a_n \phi_n( x - n )
    \,,
  \end{equation}
  then Theorem~\ref{t:existence} will still hold for this~$\omega_0$.
  All that is needed in the proof of Theorem~\ref{t:existence} is that the velocity field associated with the vortex blobs~$\set{\phi_n}$ satisfies a uniform log-Lipschitz decay bound in the form of~\eqref{e:vLogLipschitz} below.
  The proof of~\eqref{e:vLogLipschitz} provided in Appendix~\ref{s:vLogLip} for a vortex patch on the unit square can readily be adapted to this setting.
\end{remark}

\begin{example}
  This example illustrates the limitations of the well-posedness theory currently available for non-decaying weak solutions.
  Instead of using unit square vortex patches to construct~$\omega_0$, suppose we used infinite strips.
  Explicitly, consider the initial vorticity defined by
  \begin{equation}
    \omega_0 \defeq \sum_{n \in \Z} a_n \one_{S_n} \,,
  \end{equation}
  where $S_n = [n, n+1) \times \R$ is a strip,  and~$\set{a_n}$ is a family of independent, uniformly bounded, mean~$0$, variance~$1$ random variables.
  In this case we can construct an incompressible~$u_0$ with~$\crl u_0 = \omega_0$ explicitly, using shear flows.
  Indeed, one can readily verify that that the desired velocity field~$u_0$ is given by
  \begin{equation}\label{e:stripShear}
    u_0(x) =
      \begin{pmatrix}
	0\\
	a_n (x_1 - n) + \sum_{k = 0}^{n-1} a_k
      \end{pmatrix}
      \quad\text{for } x_1 \in [n, n+1)
      \,.
  \end{equation}
  This is of course a stationary solution to the Euler equations~\eqref{e:euler}.

  Since~$\sum_1^{n-1} a_k$ is a random walk after~$n$ steps, the above formula for $u_0$ shows that~$\abs{u_0(x)} = O(\abs{x_1}^{1/2})$ as~$\abs{x_1} \to \infty$.
  This falls outside the scope of the available existence theory, and as a result, we do not know if the stationary solution~$u_0$ is the unique weak solution to~\eqref{e:euler} in~$L^2_{1/2}$.
  More interestingly, if we used a deterministic bounded perturbation of the strips~$\set{S_n}$, we would still expect the~$\abs{x}^{1/2}$ growth as~$\abs{x} \to \infty$, but we would not have the luxury of the explicit stationary solution~\eqref{e:stripShear}.
  In this case even existence of a solution is not guaranteed by available results.
\iffalse
For $a,b$ real constants, a vector field of the form $(0, a x_1 + b)$ is divergence-free and has curl constant equal $a$. We can construct a continuous random velocity field $u_0$, divergence-free and with curl equal to $\omega_0$ in the sense of distributions of the form $(0,\xi_n x_1 + b_n)$, with
  $$b_n = - n\xi_n + \sum_0^{n-1} \xi_i.$$
This velocity field is a linear shear in each strip  $[n, n+1) \times \R$, and it is exactly equal to $(0,\sum_1^n \xi_i)$ at $x_1 = n+1$, for each $n\in \Z$.
This means that the velocity at position $x_1 = n+1$ is distributed as a classical 1-D random walk after $n$ time steps, and therefore its expectation is
$\mathcal{O}(\sqrt{n})$, which in turn implies that the velocity grows like $\mathcal{O}(|x|^{1/2})$ for large $|x|$, almost certainly.

This specific initial data is a stationary weak solution of the Euler equations, since $u_0$ is vertical and $\nabla \omega_0$ is horizontal, so existence of a weak solution for positive time is not an issue, but its growth falls precisely outside the known conditions for existence. For example, another velocity obtained by doing a deterministic, bounded and compactly supported perturbation of $\omega_0$ falls outside the scope of the currently available theory for well-posedness.
\fi
\end{example}

\section{Proof of Theorem \ref{t:existence}} \label{proof}

As mentioned earlier, to prove existence of solutions to the Euler equations~\eqref{e:euler} with bounded initial vorticity, we need to ensure that the corresponding initial velocity does not grow too rapidly.
%For the initial vorticity~\eqref{e:omega0def}, however, we may get unlucky in our choice of the random signs~$a_n$ resulting in a velocity field that grows too rapidly at infinity.
%This could happen, for instance, if the signs of~$a_n$ align on a sector.
We will produce an  initial velocity field $u_0$ such that, with probability~$1$, $u_0 \in L^2_\alpha$ for every~$\alpha > 0$, and then use the Cobb--Koch theorem~\cite{CobbKoch24} to obtain Theorem~\ref{t:existence}.
%will (almost surely) imply global well-posedness of the Euler equations~\eqref{e:euler}.
%(A precise growth bound is stated in Proposition~\ref{p:u0exist}, below.)

At first sight, constructing~$u_0$ might seem easy: simply use the Biot--Savart law.
That is let~$K$ be the Biot--Savart kernel given by
\begin{equation}\label{e:BiotSavart}
  K(z) = \frac{z^\perp}{2 \pi \abs{z}^2}
  \quad
  \text{where } z = (z_1, z_2) \in \R^2\,,
  \quad\text{and}\quad
  z^\perp \defeq (-z_2, z_1)
  \,,
\end{equation}
and define
\begin{equation}\label{e:u0BiotSavartBad}
  u_0 = K * \omega_0\,.
\end{equation}
When~$\omega_0 \in L^1(\R^2) \cap L^\infty(\R^2)$, we know~$\curl u_0 = \omega_0$.
(Here we used the notation~$\curl u_0 = \crl u_0 = \partial_1 u_2 - \partial_2 u_1$.)
Unfortunately, this does not work in our situation as we only know~$\omega_0$ is bounded, not integrable.
In this case~$K$ does not decay rapidly enough for the convolution in~\eqref{e:u0BiotSavartBad} to be convergent, and the construction of~$u_0$ requires a little finesse.
This is is the content of the next proposition.
%and postpone the proof to Section~\ref{s:u0exist}, below.

\begin{proposition}\label{p:u0exist}
  There exists a
  %locally log-Lipschitz
  % 2025-07-11 GI: The vector field is NOT locally log-Lipschitz!!
  random velocity field~$u_0$ such that, with probability~$1$, $u_0 \in L^2_\loc(\R^2)$, and the distributional derivative of~$u_0$ satisfies
  \begin{equation}\label{e:curlUEqOmega}
    \dv u_0 = 0
    \quad\text{and}\quad
    \curl u_0 = \omega_0
    \,,
    \quad\text{in } H^{-1}_\loc(\R^2)
    \,.
  \end{equation}
  Moreover, for every~$p \in [1, \infty)$ we have
  \begin{align}
    \MoveEqLeft
    \paren[\big]{ \E \abs{u_0(x) - u_0(y)}^p }^{1/p}
    \\
    \label{e:u0lipschitz}
      &\leq \frac{C_p \abs{x - y} (1 + \ln^{-} \abs{x - y} )  }{1 + \abs{x - y}}
	\paren[\big]{ \ln(e + \abs{x} + \abs{y}) }^{1/2}
      \,,
  \end{align}
  where~$\ln^{-}(z) = -\min \set{0, \ln z}$.
  The vector field~$u_0$ is unique up to addition of constants.
\end{proposition}

\begin{remark}[Optimality]\label{r:optimality}
  Translating~$u_0$ so that~$u_0(0) = 0$, setting~$p = 2$ and~$y = 0$ in~\eqref{e:u0lipschitz} immediately gives the square-root log growth bound
  \begin{equation}\label{e:u0Growth}
    \E \abs{u_0(x)}^2 \leq C \ln (e + \abs{x})
    \,.
  \end{equation}
  The following argument suggests that the upper bound~\eqref{e:u0Growth} may be optimal.
  For any square of side length~$n$, the probability of having~$\omega_0 \equiv 1$ on this square is~$\cramped{2^{-n^2}}$.
  Thus if we consider a region of size~$\cramped{n^2 2^{n^2}}$, it is likely we will see an~$n \times n$ square on which $\omega_0 \equiv 1$.
  The velocity in this~$n \times n$ square should grow linearly (like a rigid rotation) as we approach the center.
  This suggests that in a region of size~$\cramped{n^2 2^{n^2}}$ the velocity should grow by~$n$, which in turn implies~$\abs{u(x)}$ grows like~$\sqrt{\ln \abs{x}}$ as~$\abs{x} \to \infty$.
  This heuristic indicates that~\eqref{e:u0Growth} may be optimal.
\end{remark}

We prove Proposition~\ref{p:u0exist} in the next two subsections, and then we prove Theorem~\ref{t:existence}.

\subsection{Existence of an initial velocity}\label{s:u0exist}

%The main ingredient in the proof of Theorem~\ref{t:existence} is the construction of an initial velocity field~$u_0$ which does not grow too rapidly at infinity.
%This is the content of Proposition~\ref{p:u0exist}, and will be proved in this section.
As observed earlier, one can not directly construct~$u_0$ using the Biot--Savart law~\eqref{e:BiotSavart}.
Alternately, one can attempt to construct~$u_0$ as a linear combination of the velocity field associated to the vortex patch on the unit square~$[0, 1]^2$.
That is, we use the Biot--Savart law to define
\begin{equation}\label{e:v}
  v(x)
    = K * \one_{Q_0}(x)
    = \frac{1}{2\pi} \int_{Q_0} \frac{(x - y)^\perp}{\abs{x - y}^2} \, dy
  \,.
\end{equation}
Since~$Q_0$ is compact, \eqref{e:v} is well defined and gives a bounded incompressible velocity field~$v$ for which~$\curl v = \one_{Q_0}$.
By linearity and translation invariance, one would expect
\begin{equation}\label{e:u0first}
  u_0(x) = \sum_{n \in \Z^2} a_n v(x + n)
\end{equation}
to be the desired velocity field for which~\eqref{e:curlUEqOmega} holds.
However, it is easy to see that
\begin{equation}
  \abs{v(x+n)}
    \approx \frac{1}{\abs{x + n}}
    \,,
\end{equation}
which does not decay fast enough to ensure convergence of the sum in~\eqref{e:u0first}.

We overcome this obstacle through a renormalization: we subtract a constant vector field from each term in~\eqref{e:u0first} to obtain faster decay of each term, without affecting the validity of~\eqref{e:curlUEqOmega}.
Explicitly, define
\begin{equation}\label{e:u0}
  u_0(x) \defeq
    \sum_{n \in \Z^2}
      a_n \paren{
	v(x+n) - v(n)
      }
    \,.
\end{equation}

Using the Biot--Savart law it is easy to show (see Lemma~\ref{l:vLogLipschitz}, below) that
\begin{equation}
  \abs{v(x+n) - v(n)}
    \approx \frac{C(x)}{1 + \abs{n}^2}
    \,.
\end{equation}
Hence the sum in~\eqref{e:u0} diverges absolutely at every point, which may appear discouraging.
% which may appear discouraging as this 
% which diverges when summed over all~$n \in \Z^2$.
% It is easy to see that 
% %and~$\abs{v(x+n) - v(x)}$ does not decay fast enough to ensure the
% the series on the right hand side of~\eqref{e:u0} is not pointwise absolutely convergent.
% %Thus the series in~\eqref{e:u0} does not converge absolutely.
The fact that~$a_n$'s are chosen independently, however, provides the additional cancellation needed to ensure that the right hand side of~\eqref{e:u0} converges on average, in a sense that will be made precise shortly.
%This will be used to prove Proposition~\ref{p:u0exist}.
To carry out the details, we state a lemma concerning the decay of the velocity field~$v$.

\begin{lemma}\label{l:vLogLipschitz}
  For every~$x, y \in \R^2$ the velocity field~$v$ defined by~\eqref{e:v} satisfies the log-Lipschitz bound
  \begin{equation}\label{e:vLogLipschitz}
    \abs{v(x) - v(y)} \leq \frac{C \abs{x - y} (1 + L(x, y)) }{(1 + \abs{x})(1 + \abs{y})}
    \,.
  \end{equation}
  Here~$L(x, y)$ is defined by
  \begin{equation}\label{e:Ldef}
    L(x, y) \defeq \one_{ \set{x, y \in 3Q_0}} \ln^{-} \abs{x - y}
  \end{equation}
  where~$3Q_0$ denotes the concentric triple of~$Q_0$.
\end{lemma}

The proof of Lemma~\ref{l:vLogLipschitz} is relatively standard and is presented in Appendix~\ref{s:vLogLip} for completeness.
Notice~\eqref{e:u0} and independence of the random variables~$\set{a_n}$ immediately implies
\begin{equation}
  \E \abs{u_0(x)}^2
    = \sum_{n \in \Z^2} \abs{v(x + n) - v(n)}^2
    \,.
\end{equation}
Lemma~\ref{l:vLogLipschitz} implies that~$\abs{v(x + n) - v(n)}^2$ decays like~$1 / \abs{n}^4$, which is summable in~$\Z^2$.
Thus for every~$x \in \R^2$, we have $\E \abs{u_0(x)}^2 < \infty$, and hence the sum on the right of~\eqref{e:u0} converges (conditionally) almost surely.
The event on which this sum converges, however, may depend on~$x$.
Since there are uncountably many~$x \in \R^2$, some care has to be taken to ensure~$u_0$ is defined in a meaningful way almost surely.

The standard approach is to construct a continuous modification \`a la Kolmogorov--Censtov~\cite{KaratzasShreve91}.
While this can likely be accomplished in our situation, it is more convenient to find~$u_0$ in~$L^2_\loc(\R^2)$ as opposed to pointwise.
Of course, if~$u_0 \in L^2_\loc(\R^2)$ almost surely, standard elliptic regularity and~\eqref{e:curlUEqOmega} will imply~$u_0$ is also continuous almost surely.

To construct~$u_0$, let~$(\Omega, \mathcal F, \P)$ be the probability space on which the random variables~$\{a_n\}$ are defined.
Let~$\mathcal L^2$ be the Fr\'echet space of jointly measurable functions on~$\Omega \times \R^2$ defined by the countable family of semi-norms
%~$\set{ \norm{\cdot}_{L^2(\Omega \times B_N)} \st N \in \N}$ where
\begin{equation}
  \norm{u}_{L^2(\Omega \times B_N)}^2
    \defeq \E \norm{u}_{L^2(B_N)}^2
      = \E \int_{\abs{x} < N} \abs{u(x)}^2 \, dx
      \,,
\end{equation}
for~$N \in \N$.

\begin{lemma}\label{l:l2conv}
  The series in~\eqref{e:u0} converges in~$\mathcal L^2$.
\end{lemma}
\begin{proof}
  We first show that for every~$N \in \N$, the partial sums of the series in~\eqref{e:u0} are Cauchy in~$L^2(\Omega \times B_N)$.
  For this for every~$n \in \Z^2$ define
  \begin{equation}\label{e:vnDef}
    v_n(x) = v(x + n) - v(n)\,.
  \end{equation}
  For any~$M_1,  M_2 \in \N$ with~$3 \leq M_1 < M_2$, we note
  \begin{equation}\label{e:sumAnVn}
    \norm[\bigg]{\sum_{M_1 \leq \abs{n} \leq M_2} a_n v_n }^2_{L^2(\Omega \times B_N)}
      = \int_{B_N} \E \abs[\bigg]{\sum_{M_1 \leq \abs{n} \leq M_2} a_n v_n }^2 \, dx
  \end{equation}
  % For every~$x \in \R^d$, note
  % \begin{equation}\label{e:Eu021}
  %   \E \abs{u_0(x)}^2
  %       = \sum_{ m, n \in \Z^2 }
  %           \E \paren{a_m a_n}
  %             \paren{ v(x+m) - v(m) }^2
  %             \paren{ v(x+n) - v(n) }^2
  %       \,.
  % \end{equation}
  Since the~$a_n$'s are independent with mean~$0$ and variance~$1$ we know
  \begin{equation}
    \E (a_m a_n)
      = \begin{cases}
	  0 & m \neq n\,,\\
	  1 & m = n
	  \,.
      \end{cases}
  \end{equation}
  Using this in~\eqref{e:sumAnVn} implies
  \begin{align}
    \norm[\bigg]{\sum_{M_1 \leq \abs{n} \leq M_2} a_n v_n }^2_{L^2(\Omega \times B_N)}
      &= \int_{B_N} \sum_{M_1 \leq \abs{n} \leq M_2} \abs{v_n}^2 \, dx
    \\
      &\leq
	C \int_{B_N} \sum_{M_1 \leq \abs{n} \leq M_2} 
	  \frac{\abs{x}^2 }{ (1 + \abs{x + n} )^2 ( 1 + \abs{n} )^2 }
      \, dx
      \,.
  \end{align}
  The last inequality above followed from~\eqref{e:vLogLipschitz} and the fact that~$L(x+n, n) = 0$ as~$n \not\in 3Q_0$.
  Consequently
  \begin{equation}\label{e:psum1}
    \norm[\bigg]{\sum_{M_1 \leq \abs{n} \leq M_2} a_n v_n }^2_{L^2(\Omega \times B_N)}
      \leq
	\sum_{M_1 \leq \abs{n} \leq M_2} 
	  \frac{C N^4}{ (1 + \abs{n} - N )^2 ( 1 + \abs{n} )^2 }
      \,.
  \end{equation}
  Since
  \begin{equation}
    \sum_{\substack{n \in \Z^2,\\ \abs{n} > N}}
      \frac{1}{ (1 + \abs{n} - N )^2 ( 1 + \abs{n} )^2 } < \infty
      \,,
  \end{equation}
  the right hand side of~\eqref{e:psum1} can be made arbitrarily small by making~$M_1$ and~$M_2$ sufficiently large.
  This shows the partial sums of the right hand side of~\eqref{e:u0} are Cauchy with respect to each of the semi-norms~$\norm{\cdot}_{L^2(\Omega \times B_N)}$.
  This immediately shows the series in~\eqref{e:u0} converges in~$\mathcal L^2$, concluding the proof.
\end{proof} 

In light of Lemma~\ref{l:l2conv}, the function~$u_0$ is a jointly measurable function on~$\Omega \times \R^2$, and hence it's slices must almost surely be in~$L^2_\loc(\R^2)$.
Combined with separability of~$H^{-1}_\loc(\R^2)'$ this will imply the identities in ~\eqref{e:curlUEqOmega}.

\begin{lemma}\label{l:curlUEqOmega}
  There exists a null set~$\mathcal N \subseteq \Omega$ such that on~$\mathcal N^c$ the function~$u_0$ satisfies~\eqref{e:curlUEqOmega}.
\end{lemma}
\begin{proof}
  Let~$v_n$ be the functions defined in~\eqref{e:vnDef}.
  We first claim that for any (deterministic) test function~$\varphi \in C^\infty_c(\R^2)$, there exists a null set~$\mathcal N_\varphi$ (depending on~$\varphi$) such that on~$\mathcal N_\varphi^c$ we have
  \begin{equation}\label{e:intU}
    \sum_{n \in \Z^2} a_n \int_{\R^2} v_n \varphi \, dx
      = \int_{\R^2} u_0 \varphi \, dx
    \,.
  \end{equation}

  To prove this choose~$R > 0$ to be large enough so that~$\supp(\varphi) \subseteq B_R$.
  By Fubini's theorem we know the integral $\int_{\R^2} u_0 \varphi \, dx$ is defined almost surely, and is a square integrable random variable.
  Next, for any~$N \in \N$, $N>3$,  we note
  \begin{align}
    \MoveEqLeft
    \E\paren[\Big]{
      \int_{\R^2} u_0 \varphi \, dx - \sum_{\abs{n} < N} \int_{\R^2} a_n v_n \varphi \, dx
    }^2
      = \E\paren[\Big]{
	\int_{\R^2} \sum_{\abs{n} \geq N} a_n v_n \varphi \, dx
      }^2
    \\
      &\leq \abs{B_R}
	\int_{B_R} \E \paren[\bigg]{\sum_{\abs{n} \geq N} a_n v_n }^2 \varphi^2 \, dx
    \\
      &\leq C \abs{B_R} \norm{\varphi}_{L^\infty}^2
	\sum_{\abs{n} \geq N} \frac{R^4}{(1 + \abs{n} - R)^2 (1 + \abs{n})^2}
    \xrightarrow{N \to \infty} 0
    \,.
  \end{align}
  This immediately implies the existence of a null set~$\mathcal N_\varphi$ such that~\eqref{e:intU} holds on~$\mathcal N_\varphi^c$.
  \smallskip

  Now choose a countable set of (deterministic) smooth, compactly supported functions~$\set{\varphi_k \st k \in \N}$ which are dense in the dual space~$H^{-1}_\loc(\R^2)'$.
  For each~$k \in \N$, the above argument and~\eqref{e:intU} imply that there exists a null set~$\mathcal N_k$ such that on~$\mathcal N_k^c$ we have
  \begin{align}
    - \int_{\R^2} u_0 \cdot \grad^\perp \varphi_k \, dx
      &= - \sum_{n \in \Z^2} a_n \int_{\R^2} v_n \grad^\perp \varphi_k \, dx
      = \sum_{n \in \Z^2} a_n \int_{\R^2} \curl v_n \varphi_k \, dx
    \\
    \label{e:phikNk}
      &= \int_{\R^2} \omega_0 \varphi_k \,.
  \end{align}

  Using Lemma~\ref{l:l2conv} find a null set~$\mathcal N'$ such that $u_0 \in L^2_\loc(\R^2)$ on the complement.
  Let~$\mathcal N  = \cup_{k \geq 0} \mathcal N_k \cup \mathcal N'$.
  On~$\mathcal N^c$, we note $u_0 \in L^2_\loc(\R^2)$ and hence~$\curl u_0 \in H^{-1}_\loc(\R^2)$.
  Moreover, \eqref{e:phikNk} shows that for every~$k \in \N$ we have
  \begin{equation}
    \int_{\R^2} \curl u_0 \varphi_k \, dx
      = \int_{\R^2} \omega_0 \varphi_k \, dx
      \quad\text{on } \mathcal N^c
      \,.
  \end{equation}
  Since~$\curl u_0 \in H^{-1}_\loc(\R^2)$, and~$\set{\varphi_k}$ is dense in~$H^{-1}_\loc(\R^2)'$, this implies that~$\curl u_0 = \omega_0$ in~$H^{-1}_\loc(\R^2)$.
  A similar argument shows~$\dv u_0 = 0$ in~$H^{-1}_\loc(\R^2)$, thereby concluding the proof.
\end{proof}

\subsection{Growth of the initial velocity (Proposition \ref{p:u0exist})}
We are now in a position to prove Proposition~\ref{p:u0exist}.
By Lemmas \ref{l:l2conv} and \ref{l:curlUEqOmega} we already have (almost sure) existence of a divergence free velocity field~$u_0$ with~$\curl u_0 = \omega_0$.
The main point of Proposition~\ref{p:u0exist}, however, is to ensure the growth bound~\eqref{e:u0lipschitz}, which will, in turn, be the key to establishing an almost sure $L^2_\alpha$-norm estimate for $u_0$, needed to apply the Cobb-Koch well-posedness theorem.

\begin{proof}[Proof of Proposition~\ref{p:u0exist}]
  Lemma~\ref{l:curlUEqOmega} shows almost sure existence of an incompressible velocity field~$u_0 \in L^2_\loc(\R^2)$ that satisfies~\eqref{e:curlUEqOmega}.
  It remains to prove~\eqref{e:u0lipschitz}.
  Since the sum in~\eqref{e:u0} is convergent, we see
  \begin{equation}
    u_0(x) - u_0(y)
      = \sum_{n \in \Z^2}
      a_n \paren{
	v(x+n) - v(y+n)
      }
    \,.
  \end{equation}
  Using Khintchine's inequality (see for instance Section 2.6 in~\cite{Vershynin18}), this implies for every~$p \in [1, \infty)$ there exists a constant~$C_p$ such that
  \begin{align}
    \paren[\big]{ \E \abs{u_0(x)  - u_0(y)}^p  }^{2/p}
      &\leq C_p \sum_{n \in \Z^2}
	\abs{ v(x+n) - v(y + n) }^2
    \\
    \label{e:Eu0xMinusU0xPrime}
       &\overset{\mathclap{\eqref{e:vLogLipschitz}}}{\leq}
	C_p \abs{x - y}^2 \paren[\big]{
	  1 + \ln^{-} \abs{x - y}
	}^2
	  S(x, y)
	  \,,
  \end{align}
  where
  \begin{equation}
    S(x, y) \defeq \sum_{n \in \Z^2} \frac{1}{ (1 + \abs{n - x})^2 ( 1 + \abs{n - y} )^2 }
  \,.
  \end{equation}

  We will now bound the sum~$S(x, y)$ by partitioning~$\Z^2$ into five regions.
  Without loss of generality we assume~$\abs{x} \leq \abs{y}$, and divide the analysis into two cases.

  \restartcases
  \case[$\abs{x - y} \geq \abs{y} / 2$]
  Define the sets~$A_1$, \dots, $A_5$ by
  \begin{subequations}
  \begin{align}
    \label{e:A1def}\noeqref{e:A1def}
    A_1 &= \set[\Big]{ n \in \Z^2 \st \abs{n} \leq \frac{\abs{x}}{2} }
  \\
    \label{e:A2def}\noeqref{e:A2def}
    A_2 &= \set[\Big]{ n \in \Z^2 \st \abs{n} \geq 2 \abs{y} }
  \\
    \label{e:A3def}\noeqref{e:A3def}
    A_3 &= \set[\Big]{ n \in \Z^2 - A_1 - A_2 \st \abs{n - x} \leq \frac{\abs{x}}{4} }
  \\
    \label{e:A4def}\noeqref{e:A4def}
    A_4 &= \set[\Big]{ n \in \Z^2 - A_1 - A_2 - A_3 \st \abs{n - y} \leq \frac{\abs{y}}{4} }
  \\
    \label{e:A5def}\noeqref{e:A5def}
    A_5 &= \Z^2 - \bigcup_{i = 1}^4 A_i
    \,.
  \end{align}
  \end{subequations}

  Clearly for all~$n \in A_1$ we have
  \begin{equation}\label{e:ninA1}
    \abs{n - x} \geq \frac{\abs{x}}{2}
    \quad\text{and}\quad
    \abs{n - y} \geq \frac{\abs{y}}{2}
  \end{equation}
  and so
  \begin{equation}\label{e:sumA1}
    \sum_{n \in A_1}
      \frac{1}{ (1 + \abs{n - x})^2 ( 1 + \abs{n - y} )^2 }
      \leq \frac{C \abs{A_1}}
	{\paren{ 1 + \abs{x}^2 } \paren{ 1 + \abs{y}^2 }}
	\leq \frac{C}{1 + \abs{y}^2}
      \,.
  \end{equation}

  Next for all~$n \in A_2$ we have
  \begin{equation}\label{e:ninA2}
    \abs{n - x} \geq \frac{\abs{n}}{2}
    \quad\text{and}\quad
    \abs{n - y} \geq \frac{\abs{n}}{2}
  \end{equation}
  and so
  \begin{align}
    \sum_{n \in A_2}
      \frac{1}{ (1 + \abs{n - x})^2 ( 1 + \abs{n - y} )^2 }
      &\leq \sum_{\abs{n} \geq 2 \abs{y}}
	\frac{C}{1 + \abs{n}^4}
	\leq C \int_{\abs{z} \geq 2 \abs{y}} \frac{dz}{1 + \abs{z}^4}
    \\
      \label{e:sumA2}\noeqref{e:sumA2}
	&\leq \frac{C}{1 + \abs{y}^2}
      \,.
  \end{align}

  Next for~$n \in A_3$ we note
  \begin{equation}\label{e:ninA3}
    \abs{n - y}
      \geq \abs{y - x} - \abs{n - x}
      \geq \frac{\abs{y}}{4}
    \,,
  \end{equation}
  and so
  \begin{align}
    \sum_{n \in A_3}
      \frac{1}{ (1 + \abs{n - x})^2 ( 1 + \abs{n - y} )^2 }
      &\leq
	\frac{16}{1 + \abs{y}^2} \int_{\abs{z - x} \leq \frac{\abs{x}}{4} } \frac{dz}{1 + \abs{z - x}^2}
    \\
      \label{e:sumA3}\noeqref{e:sumA3}
      &\leq
	\frac{C \ln( e + \abs{x}) }{1 + \abs{y}^2}
    \,.
  \end{align}
  The sum in~$A_4$ is handled similarly to the sum in~$A_3$ and gives
  %\GI[2025-07-11]{Not a typo, since we have~$\abs{n - x} \geq \abs{y} / 4$ in this case.}%
  \begin{equation}\label{e:sumA4}\noeqref{e:sumA4}
    \sum_{n \in A_4}
      \frac{1}{ (1 + \abs{n - x})^2 ( 1 + \abs{n - y} )^2 }
    \leq \frac{C \ln( e + \abs{y}) }{1 + \abs{y}^2}
    \,.
  \end{equation}

  Finally, for~$n \in A_5$, we note
  \begin{equation}\label{e:ninA5}
    \abs{n - x} \geq \frac{\abs{x}}{4}
    \quad\text{and}\quad
    \abs{n - y} \geq \frac{\abs{y}}{4}
    \,,
  \end{equation}
  and hence
  \begin{equation}\label{e:sumA5}
    \sum_{n \in A_5}
      \frac{1}{ (1 + \abs{n - x})^2 ( 1 + \abs{n - y} )^2 }
    \leq
      \frac{C \abs{A_5}}{(1 + \abs{y})^4}
    \leq
      \frac{C(\abs{y}^2 - \abs{x}^2)}{(1 + \abs{y})^4}
      \,.
  \end{equation}
  \medskip

  Combining~\eqref{e:sumA1}, \eqref{e:sumA2}, \eqref{e:sumA3}, \eqref{e:sumA4} and~\eqref{e:sumA5} shows
  \begin{equation}\label{e:Scase1}
    S(x, y)
      \leq \frac{C \ln(e + \abs{y})}{1 + \abs{y}^2}
      \leq \frac{C \ln(e + \abs{y})}{1 + \abs{x - y}^2}
      \,,
  \end{equation}
  where the last inequality follows from the fact that~$\abs{x - y} \leq 2 \abs{y}$.
  \medskip

  \case[$\abs{x - y} < \abs{y} / 2$]
  The analysis in this case is very similar to the previous case, except we need to change the sets~$A_3$, $A_4$ and~$A_5$.
  Let~$A_1$, $A_2$ be as in~\eqref{e:A1def} and~\eqref{e:A2def} respectively, and define
  \begin{subequations}
  \begin{align}
    \label{e:A3defP}\noeqref{e:A3defP}
    A_3' &= \set[\Big]{ n \in \Z^2 - A_1 - A_2 \st \abs{n - x} \leq \frac{\abs{x - y}}{2} }
  \\
    \label{e:A4defP}\noeqref{e:A4defP}
    A_4' &= \set[\Big]{ n \in \Z^2 - A_1 - A_2 - A_3 \st \abs{n - y} \leq \frac{\abs{x-y}}{2} }
  \\
    \label{e:A5defP}\noeqref{e:A5defP}
    A_5' &= \Z^2 - A_1 - A_2 - A_3' - A_4'
    \,.
  \end{align}
  \end{subequations}

  For~$n \in A_3'$ we note
  \begin{equation}\label{e:ninA3P}
    \abs{n - y}
      \geq \frac{\abs{x - y}}{2}
    \,,
  \end{equation}
  and so
  \begin{align}
    \sum_{n \in A_3'}
      \frac{1}{ (1 + \abs{n - x})^2 ( 1 + \abs{n - y} )^2 }
      &\leq
	\frac{4}{1 + \abs{x - y}^2} \int_{\abs{z - x} \leq \frac{\abs{x - y}}{2} } \frac{dz}{(1 + \abs{z - x})^2}
    \\
      \label{e:sumA3P}\noeqref{e:sumA3P}
      &\leq
	\frac{C \ln( e + \abs{x - y}) }{(1 + \abs{x - y})^2}
    \,.
  \end{align}
  The bound on~$A_4'$ is similar and gives the exact same bound as~\eqref{e:sumA3P}.

  Finally, for~$n \in A_5'$ we note that both
  \begin{equation}\label{e:ninA5p}
    \abs{n - x} \geq \frac{\abs{x -y}}{2}
    \quad\text{and}\quad
    \abs{n - y} \geq \frac{\abs{x -y}}{2}
    \,,
  \end{equation}
  and so
  \begin{equation}\label{e:sumA5P}
    \sum_{n \in A_5'}
      \frac{1}{ (1 + \abs{n - x})^2 ( 1 + \abs{n - y} )^2 }
    \leq
      \frac{16 \abs{A_5}}{(1 + \abs{x - y})^4}
    \leq
      \frac{C(\abs{y}^2 - \abs{x}^2)}{(1 + \abs{x - y})^4}
      \,.
  \end{equation}

  Combining~\eqref{e:sumA3P}, and~\eqref{e:sumA5P} shows
  \begin{equation}\label{e:Scase2}
    S(x, y)
      \leq \frac{C \ln( e + \abs{ x - y} )}{1 + \abs{x - y}^2}
      \leq \frac{C \ln( e + \abs{y} )}{1 + \abs{x - y}^2}
    \,,
  \end{equation}
  where the last inequality followed from the assumption~$\abs{x - y} < \abs{y} / 2$.
  \medskip

  Using~\eqref{e:Scase1} and~\eqref{e:Scase2} in~\eqref{e:Eu0xMinusU0xPrime} yields~\eqref{e:u0lipschitz}, completing the proof.
\end{proof}
\subsection{Proof of global well-posedness (Theorem~\ref{t:existence})}

We now prove Theorem~\ref{t:existence}.
As mentioned earlier, for global well-posedness of the Euler equations in~$\R^2$, one needs both bounded initial vorticity, and a growth condition on the initial velocity.
We emphasize that, while~\eqref{e:u0Growth} provides a bound on the growth of~$\E \abs{u_0}^2$, it does not guarantee an almost sure bound on the growth of~$u_0$.
That is, the bound~\eqref{e:u0Growth} \emph{does not} directly imply the existence of a (random) constant~$\widetilde{C}$, with~$\P(\widetilde{C} < \infty) = 1$, such that, for every~$x \in \R^d$, we have the \emph{pointwise} bound
  \begin{equation}\label{e:u0Pointwise}
  \abs{u_0(x)}^2 \leq \widetilde{C} \ln( e + \abs{x})
  \,.
\end{equation}
%As a result, Theorem~\ref{t:existence} does not directly follow from results such as~\cite{Cozzi15,CozziKelliher19} which require on pointwise bounds on the growth of the initial data.

The bound~\eqref{e:u0Growth}, however, does imply growth bounds that are averaged in space.
Indeed, Proposition~\ref{p:u0exist} will imply~$u_0 \in L^2_\alpha$ almost surely, after which the Cobb--Koch well-posedness  result~\cite{CobbKoch24} will yield Theorem~\ref{t:existence}.

\begin{proof}[Proof of Theorem~\ref{t:existence}]
  We recall the result in~\cite[Theorem 1]{CobbKoch24}, establishing the  existence and uniqueness of a global weak solution of the Euler equations provided the initial data belongs to~$L^2_\alpha$ for some~$\alpha \in [0, 1/2)$, and the initial vorticity is bounded.
  In our case, the boundedness of vorticity is immediate.
  Indeed, the random variables~$a_n$ are all uniformly bounded by assumption, and so
  \begin{equation}
    \norm{\omega_0}_{L^\infty} \leq \max_n \max \abs{a_n} < \infty
    \,.
  \end{equation}
  %and so we only need to ensure the initial velocity field~$u_0$ is in~$L^2_\alpha$.
  We will now show that for every~$\alpha > 0$, Proposition~\ref{p:u0exist} implies that $u_0 \in L^2_\alpha$ almost surely.

  To see this, let~$R \geq 1$, and $\alpha > 0$ be arbitrary.
  Observe
  \begin{align}
    \frac{1}{R^{2 + 2 \alpha }} \int_{\abs{x} \leq R} \abs{u_0}^2 \, dx
      &\leq C
	\int_{\abs{x} \leq R} \frac{ \abs{u_0(x)}^2 }{1 + \abs{x}^{2 + 2\alpha}} \, dx
      \leq C
	\int_{\R^2} \frac{ \abs{ u_0(x) }^2 }{1 + \abs{x}^{2 + 2\alpha}} \, dx
      \,.
  \end{align}
  Thus
  \begin{align}
    \E \norm{u_0}_{L^2_\alpha}^2
      &= \E \sup_{R \geq 1}
	\frac{1}{R^{2 + 2 \alpha }} \int_{\abs{x} \leq R} \abs{u_0}^2 \, dx
      \leq C
	\E \int_{\R^2} \frac{ \abs{ u_0(x) }^2 }{1 + \abs{x}^{2 + 2\alpha}} \, dx
    \\
      &\overset{\mathclap{\eqref{e:u0Growth}}}{\leq} ~C
	\int_{\R^2} \frac{\ln( e + \abs{x})  }{1 + \abs{x}^{2 + 2\alpha}} \, dx
      < \infty
      \,.
  \end{align}
  This implies~$\E \norm{u_0}_{L^2_\alpha} < \infty$ and so~$\norm{u_0}_{L^2_\alpha}$ is finite almost surely.
  Thus Theorem~1 in~\cite{CobbKoch24} applies to almost every realization of~$u_0$, and this implies global well-posedness as desired.
\end{proof}

\section{Decorrelation and Scaling} \label{decorscaling}

\subsection{Decorrelation}

In our setting we note that the initial vorticity has a distinguished length scale: points that are more than a distance of~$\sqrt{2}$ apart have independent vorticities.
The velocity, however, is obtained in a linear but non-local manner from the initial vorticity.
Moreover, the Biot--Savart kernel decays slowly and so the influence of each vortex patch is felt on a large region.
As a result, the velocity experienced by points that are large distances apart may be highly correlated.
We now show that the covariance of the velocity fields at two different points grows logarithmically with the distance.
This is the same order of magnitude as~\eqref{e:u0Growth}, suggesting that the velocity may not decorrelate even over large scales.
An interesting question, that goes beyond the scope of the present work, is to find at positive times~$ t > 0$ a distinguished length scale of the vorticity, or the decorrelation of the velocity.

\begin{proposition}[Decorrelation bounds]
  For every~$x, y \in \R^2$ we have
  \begin{equation}
    \abs{\cov( u_0(x),  u_0(y))}
      \leq C \ln(e + \abs{x - y})
      \,,
      \quad\text{provided}\quad
      \abs{x - y} \geq \max \set[\Big]{ \frac{\abs{x}}{2}, \frac{\abs{y}}{2} }
      \,.
  \end{equation}
\end{proposition}
\begin{proof}
  Using~\eqref{e:u0}, we note~$\E u_0(x) = \E u_0(y) = 0$ and hence
  \begin{align}
    \abs{\cov( u_0(x),  u_0(y))}
      &= \E u_0(x) \cdot u_0(y)
      = \sum_{n \in \Z^2} ( v(x+n) - v(n) ) \cdot (v(y+n) - v(n))
    \\
    \label{e:covu0xu0y}
      &\overset{\mathclap{\eqref{e:vLogLipschitz}}}{\leq}
	C \abs{x} \abs{y} (1 + \ln^{-} \abs{x})(1 + \ln^{-} \abs{y})
	\sum_{n \in \Z^2} F_n(x, y)
    \,,
  \end{align}
  where
  \begin{equation}
    F_n(x, y)
      \defeq \frac{1}{(1 + \abs{x - n})(1 + \abs{y - n}) (1 + \abs{n})^2}
    \,.
  \end{equation}
  We now bound~$\sum F_n(x, y)$ by dividing~$\Z^2$ into five regions.
  Without loss of generality we assume~$\abs{x} \leq \abs{y}$.
  %and divide the analysis into two cases.

  %\restartcases
  %\case[$\abs{x - y} \geq \abs{y} / 2$]
  Let~$A_1$, \dots, $A_5$ be as in~\eqref{e:A1def}--\eqref{e:A5def}.
  For~$n \in A_1$ we use~\eqref{e:ninA1} and obtain
  \begin{equation}\label{e:sumA1cov}
    \sum_{n \in A_1}
      F_n(x, y)
      \leq \frac{C}{\paren{ 1 + \abs{x} } \paren{ 1 + \abs{y} }}
	\int_{\abs{z} \leq \frac{\abs{x}}{2}} \frac{dz}{(1 + \abs{z})^2}
	\leq \frac{C \ln( e + \abs{x} )}{(1 + \abs{x})(1 + \abs{y})}
      \,.
  \end{equation}

  Next for~$n \in A_2$ we use~\eqref{e:ninA2} and obtain
  \begin{equation}\label{e:sumA2cov}
    \sum_{n \in A_2}
      F_n(x, y)
      %\leq C \sum_{\abs{n} \geq 2\abs{y}} \frac{1}{(1 + \abs{n})^4}
      \leq C \int_{\abs{z} \geq 2\abs{y}} \frac{dz}{(1 + \abs{z})^4}
      \leq \frac{C}{(1 + \abs{y})^2}
      \,.
  \end{equation}

  For~$n \in A_3$ we use~\eqref{e:ninA3} and the fact that~$\abs{n} \geq \abs{x}/2$ to obtain
  \begin{equation}\label{e:sumA3cov}
    \sum_{n \in A_3} F_n(x, y)
      \leq
	\frac{C}{(1 + \abs{x})^2 (1 + \abs{y})}
	\int_{\abs{z - x} \leq \frac{\abs{x}}{4}}
	  \frac{dz}{1 + \abs{z}}
      \leq
	\frac{C  }{(1 + \abs{x}) (1 + \abs{y})}
      \,.
  \end{equation}

  The bound for the sum in~$A_4$ is similar and gives
  \begin{equation}
    \sum_{n \in A_4} F_n(x, y)
      \leq \frac{C }{(1 + \abs{x}) (1 + \abs{y})}
      \,.
  \end{equation}

  Finally for~$A_5$ we use~\eqref{e:ninA5} to obtain
  \begin{align}
    \sum_{n \in A_5} F_n(x, y)
      &\leq \frac{C}{(1 + \abs{x})(1 + \abs{y})}
	\int_{\frac{\abs{x}}{2} \leq \abs{z} \leq 2 \abs{y}}
	  \frac{dz}{(1 + \abs{z})^2}
      \\
      &\leq \frac{C \paren[\big]{ \ln( 1 + 2\abs{y}) - \ln(1 + \abs{x}/2) }}{(1 + \abs{x})(1 + \abs{y})}
      \leq \frac{C \ln (e + \abs{y})}{(1 + \abs{x})(1 + \abs{y})}
      \,.
  \end{align}

  Combining the above, and using the fact that~$\abs{y} \leq 2 \abs{x - y}$, we obtain
  \begin{equation}\label{e:covCase1}
    \sum_{n \in \Z^2} F_n(x, y)
      \leq \frac{C \ln( e + \abs{x - y}) }{(1 + \abs{x})(1 + \abs{y})}
      \,.
  \end{equation}
  Using this in~\eqref{e:covu0xu0y} concludes the proof.
  \iffalse
  \case[$\abs{x - y} < \abs{y} / 2$]
  As in the proof of Proposition~\ref{p:u0exist}, we partition~$\Z^2$ into the sets~$A_1$, $A_2$ defined by~\eqref{e:A1def}, \eqref{e:A2def} and the sets~$A_3'$, $A_4'$, $A_5'$ defined by~\eqref{e:A3defP}--\eqref{e:A5defP}.

  For~$n \in A_3'$ we use~\eqref{e:ninA3P} and the fact that~$\abs{n} \geq \abs{x}$ to obtain
  \begin{align}
    \sum_{n \in A_3'} F_n(x, y)
      &\leq \frac{C}{(1 + \abs{x - y}) (1 + \abs{x})^2}
	\int_{\abs{z - x} \leq \frac{\abs{x - y}}{2} }
	  \frac{dz}{1 + \abs{z - x}}
    \\
      &\leq \frac{C \abs{x - y}}{( 1 + \abs{x-y})(1 + \abs{x})^2}
      \leq \frac{C }{( 1 + \abs{x-y})(1 + \abs{x})}
      \,,
  \end{align}
  where the last inequality is true because~$\abs{x - y} \leq 4 \abs{x}$.
  The bound for~$n \in A_4'$ is identical and gives the same bound as that for~$A_3'$.

  Finally, for~$n \in A_5'$ we use~\eqref{e:ninA5p} to obtain
  \begin{equation}
    \sum_{n \in A_5'} F_n(x,y)
      \leq \frac{C}{(1 + \abs{x - y})^2} \int_{\frac{\abs{x}}{2} \leq \abs{z} \leq 2 \abs{y}} \frac{dz}{(1 + \abs{z})^2}
      %\leq \frac{\paren[\big]{ \ln( 1 + 2\abs{y}) - \ln(1 + \abs{x}/2) }}{(1 + \abs{x - y})^2}
      \leq \frac{C \ln (e + \abs{y} )}{(1 + \abs{x - y})^2}
      \,.
  \end{equation}

  Combining the above we see that
  \begin{equation}
    \sum_{n \in \Z^2} F_n(x, y) \leq \textcolor{red}{TODO}
  \end{equation}
  \fi
\end{proof}

\subsection{Scaling}

The initial vorticity chosen in~\eqref{e:omega0def} has a characteristic length scale of order~$1$.
We now investigate what happens to this under rescaling.
Let~$\epsilon > 0$ define
\begin{equation}
  \omega_0^\epsilon(x) \defeq \omega_0\paren[\Big]{ \frac{x}{\epsilon} }
  \quad\text{and}\quad
  u_0^\epsilon \defeq \epsilon u_0\paren[\Big]{\frac{x}{\epsilon}}
  \,,
\end{equation}
where~$\omega_0$ is defined by~\eqref{e:omega0def}, and~$u_0$ is defined by~\eqref{e:u0}.
Using~\eqref{e:u0lipschitz} and a direct calculation we obtain
\begin{equation}
  \E \norm{u^\epsilon_0}^2_{L^2_\alpha} \leq C \epsilon^2 \abs{\ln \epsilon}
  \xrightarrow{\epsilon \to 0} 0
  \,.
\end{equation}
Moreover, using an argument similar to that in Section~\ref{s:u0exist} we can show that, almost surely, for every~$\alpha > 0$, $u^\epsilon_0$ is in the~$L^2_\alpha$-closure of~$C^\infty_c(\R^2)$.
Since~$\omega_0^\epsilon$ are uniformly bounded in~$L^\infty$ we may apply strong continuity of the solution map (given by~\cite[Theorem 1]{CobbKoch24}) to conclude that for every~$t > 0$, we have~$u^\epsilon(t) \to 0$ in~$L^2_\alpha$ as~$\epsilon \to 0$.

An interesting question, that goes beyond the scope of the present work, is to find a scaling of the initial vorticity that gives a non-trivial limit as~$\epsilon \to 0$.

We add a final remark. As we point out in the introduction, an initial vorticity of the form \eqref{e:omega0intro}, with i.i.d. $a_n$'s, is statistically periodic, in the sense that the process used to generate it is $\Z^2$-invariant. This means that the random vorticity field is statistically {\it homogeneous} at large scales.  However, the corresponding velocity law is not $\Z^2$-invariant, and the velocity is not asymptotically (statistically) homogeneous. Deterministically, periodic vorticity generates periodic velocity if and only if it has integral zero over any given period. Our vorticity has integral zero, statistically, but random imbalances on the mass of vorticity of the blobs build up when repeated and create, at large scales, the slow decorrelation decay and logarithmic growth of the velocity which we observed. Therefore, we cannot expect the weak solution obtained to retain $\Z^2$-invariance, or the vorticity to remain statistically homogeneous at large scales, unless other dynamic mechanisms come into play. This should be investigated in future research.

\iffalse %Future work / final remarks
We add a few final remarks. Our result raises several natural questions that may lead to future research, some of which we enumerate below.
\begin{enumerate}

\item Our initial data has a distinguished scale of unit length. Does this remain true for the solution? It would be interesting to show existence of a length scale that grows in time, given, for instance by the decorrelation distance of $\omega(t,\cdot)$.

\item It would be interesting to investigate the existence of solutions to the Navier-Stokes equations and the convergence of the vanishing viscosity limit.

\item Study the mixing properties of the  corresponding random flow.

\item Characterize the large time behavior of the solutions.

\item Are these solutions limits of solutions on large square boxes, or of large periodic boxes? 

\item Is it possible to do numerical approximations of these weak solutions, or, in general, of Cobb and Koch solutions, in a sensible manner?

\end{enumerate}
\fi

\appendix
\section{A log-Lipschitz bound on \texorpdfstring{$v$}{v}}\label{s:vLogLip}

It is well known that a bounded compactly supported vorticity has a log-Lipschitz velocity.
The reason we state Lemma~\ref{l:vLogLipschitz} here is because we need the fact that, away from the support of the vorticity, the velocity field is in fact Lipschitz and not only log-Lipschitz.
The proof follows standard techniques and is only presented here for completeness.

\begin{proof}[Proof of Lemma~\ref{l:vLogLipschitz}]
  Let~$K$ be the Biot--Savart kernel~\eqref{e:BiotSavart} and recall the elementary identities
 \begin{equation}\label{e:magic}
   |K(a)-K(b)| = \frac{1}{2\pi} \frac{|a-b|}{|a||b|},
   \quad\text{and}\quad
   |K(a)| = \frac{C}{|a|}\,,
 \end{equation}
 valid for any $a$, $b\in\R^2$, $a\neq b$.

 Let $x$, $y \in \R^2$, and divide the analysis into cases.

 \restartcases
  \case[$x\in \R^2 - 3Q_0$ and $y \in \R^2 - 3Q_0$]
  Let $w \in Q_0$, so that $|w| \leq \sqrt{2}$. Then, since $|x|\geq 3\sqrt{2}$ and $|y|\geq 3\sqrt{2}$, it
  follows that
  \[|x-w| \geq \frac{1}{2}(1+|x|) \text{ and } |y-w| \geq \frac{1}{2}(1+|y|).\]
  Therefore, using also \eqref{e:magic},
  \begin{align}
  |v(x)-v(y)| & = \left| \int_{Q_0} (K(x-w)-K(y-w))dw\right| \\
              & \leq \frac{1}{2\pi}\int_{Q_0} \frac{|x-y|}{|x-w||y-w|} dw \\
              & \leq \frac{2}{\pi} \frac{|x-y|}{(1+|x|)(1+|y|)}\\
              & \leq \frac{C |x - y| (1 + L(x, y)) }{(1 + |x|)(1 + |y|)},
  \end{align}
since $L(x,y) \geq 0$ .

  \case[either $x \in 3Q_0$ and $y \in \R^2 - 3Q_0$ or $x \in \R^2 - 3Q_0$ and $y \in 3Q_0$]
  Without loss of generality assume that $x \in 3Q_0$ and $y \in \R^2 - 3Q_0$. Note that, for any $w \in Q_0$, we find
  \[|y-w|\geq \frac{1}{2}(|1+|y|) \text{ and } |x-w| \leq 4 \sqrt{2}.  \]
  Then, once again using \eqref{e:magic}, we have
  \begin{align}
  |v(x)-v(y)| & \leq \frac{1}{2\pi}\int_{Q_0} \frac{|x-y|}{|x-w||y-w|} dw \\
              & \leq \frac{1}{\pi} \frac{|x-y|}{1+|y|}\int_{\{|x-w|\leq 4\sqrt{2}\}}\frac{1}{|x-w|} dw \\
              & \leq C \frac{|x-y|}{1+|y|} \leq \frac{C |x - y| (1 + L(x, y)) }{(1 + |x|)(1 + |y|)}.
  \end{align}

  \case[$x$, $y \in 3Q_0$]
  This case is classical but, for completeness' sake, we include the proof. We split $Q_0$ into $Q_{0,1} \cup Q_{0,2}$, where
  \[Q_{0,1} \defeq \left\{ w \in Q_0 \, | \, |x-y| \geq \frac{|y-w|}{2}\right\},\]
  \[Q_{0,2} \defeq \left\{ w \in Q_0 \, | \, |x-y| < \frac{|y-w|}{2}\right\}.\]

 If $w \in Q_{0,1}$ then $|y-w| \leq 2|x-y|$ and, also, $|x-w| \leq 3|x-y|$.
 If $w \in Q_{0,2}$ then $|x-w| \geq \frac{|y-w|}{2}$ and, also, $|y-w| \leq 4\sqrt{2}$.
 Therefore,
 \begin{align}
 |v(x) - v(y)| = \left| \int_{Q_0} (K(x-w)-K(y-w))dw\right|  \leq A + B
\end{align}
where
  \begin{equation}
    A \defeq
      \int_{Q_{0,1}} (|K(x-w)|+|K(y-w)|)dw
    \quad\text{and}\quad
    B \defeq  \int_{Q_{0,2}} |K(x-w)-K(y-w)|dw
    \,.
  \end{equation}

 Let us estimate the term $A$:
 \begin{align}
 A & =\int_{Q_{0,1}} (|K(x-w)|+|K(y-w)|)dw \\
 & \leq \int_{\{|x-w|\leq 3|x-y|\}} \frac{C}{|x-w|} dw +
  \int_{\{|y-w|\leq 2|x-y|\}} \frac{C}{|y-w|} dw \\
  & \leq C|x-y| \leq  \frac{C |x - y| (1 + L(x, y)) }{(1 + |x|)(1 + |y|)}.
 \end{align}
 Next we analyze $B$:
 \begin{align}
 B & \leq \int_{Q_{0,2}} |K(x-w) - K(y-w)|dw \leq \int_{Q_{0,2}} \frac{1}{2\pi}\frac{|x-y|}{|x-w||y-w|} dw \\
   & \leq C|x-y|\int_{Q_{0,2}}  \frac{1}{ |y-w|^2} dw \\
   & \leq C|x-y|\int_{\{2|x-y| \leq |y-w| \leq 4\sqrt{2} \}}  \frac{1}{ |y-w|^2} dw \\
   & = C|x-y|(\ln (2\sqrt{2}) - \ln |x-y|) \\
   & \leq C \frac{|x-y|}{1+|y|} \leq \frac{C |x - y| (1 + L(x, y)) }{(1 + |x|)(1 + |y|)}.
   \qedhere
 \end{align}

\end{proof}

\bibliographystyle{halpha-abbrv}
\bibliography{gautam-refs1,gautam-refs2,preprints}
% DO NOT EDIT THIS LINE: $Id: f603ba329f48b45853223ef4d4008d53611b0380 $
\end{document}